\documentclass[a4paper,11pt]{article}

\usepackage{latexsym}
\usepackage{amssymb}
\usepackage{amsfonts}
\usepackage{amsmath}
\usepackage{indentfirst}
\usepackage{graphicx}
\usepackage{epsfig}
\usepackage{xypic}
\usepackage[english]{babel}

\usepackage{fullpage}

\newtheorem{theorem}{Theorem}[section]
\newtheorem{prop}[theorem]{Proposition}
\newtheorem{lemma}[theorem]{Lemma}
\newtheorem{cor}[theorem]{Corollary}

\newcommand{\cvd}{\hfill $\blacksquare$\bigskip}

\date{}

\author{Giulio Cerbai\thanks{Dipartimento di Matematica e Informatica ``U.
Dini", University of Firenze, Firenze, Italy,
\tt{giulio.cerbai@unifi.it, luca.ferrari@unifi.it}}
\and Luca Ferrari$^{\dag}$}

\title{Permutation patterns in genome rearrangement problems: the reversal model\footnote{Both authors are members of the INdAM Research group GNCS; they are partially supported by INdAM - GNCS 2018 project ``Propriet\'a combinatorie e rilevamento di pattern in strutture discrete lineari e bidimensionali" and by a grant of the "Fondazione della Cassa di Risparmio di Firenze" for the project "Rilevamento di pattern: applicazioni a memorizzazione basata sul DNA, evoluzione del genoma, scelta sociale".}}

\begin{document}

\maketitle

\begin{abstract}
In the context of the genome rearrangement problem, we analyze two well known models,
namely the reversal and the prefix reversal models,
by exploiting the connection with the notion of permutation pattern.
More specifically,
for any $k$, we provide a characterization of the set of permutations having distance $\leq k$ from the identity
(which is known to be a permutation class) in terms of what we call \emph{generating peg permutations}
and we describe some properties of its basis, which allow to compute such a basis for small values of $k$.
\end{abstract}

\bigskip

\section{Introduction}

One of the major trends in bioinformatics and biomathematics is the study of the genome rearrangement problem.
Roughly speaking, given a genome, one is interested in understanding how the genome can evolve into another genome.
To give a proper formalization, several models for rearranging a genome have been introduced,
each of which defines a series of allowed elementary operations to be performed on a genome in order to obtain an adjacent one.
For several models, it is possible to define a \emph{distance} between two genomes,
by counting the minimum number of elementary operations needed to transform one genome into the other.
The investigation of the main properties of such a distance becomes then a key point in understanding the main features of the model under consideration.

\bigskip

A common formalization of any such models consists of encoding a genome using a \emph{permutation} (in linear notation)
and describing an elementary operation as a \emph{combinatorial} operation on the entries of such a permutation.
Many genome rearrangement models have been studied under this general framework.
Among them, the following ones are very well known.

\begin{itemize}
\item The \emph{reversal} model consists of a single operation, defined as follows:
a new permutation is obtained from a given one by selecting a cluster of consecutive elements and reversing it.
More formally, given $\pi =\pi_1 \pi_2 \cdots \pi_n$,
a reversal is performed by choosing $i<j<n$ and then forming the permutation
$\sigma =\pi_1 \cdots \pi_{i-1} \boxed{\pi_j \pi_{j-1}\cdots \pi_{i+1} \pi_i}\pi_{j+1}\cdots \pi_n$.
This model was introduced in \cite{WEHM}, then studied for instance in \cite{BaPe1,HP}.
\item A variant of the reversal model is the \emph{prefix reversal model},
which is a specialization of the previous one in which the reversal operation can only be performed on a prefix of the given permutation.
This is also known as \emph{pancake sorting} (see for instance \cite{GP}).
\item A very popular and studied model is the \emph{transposition model}, see \cite{BaPe2}.
Given a permutatation $\pi =\pi_1 \cdots \pi_n$, a \emph{transposition operation} consists of
taking two adjacent clusters of consecutive elements and interchanging their positions.
Formally, one has to choose indices $1\leq i<j<k\leq n+1$,
then form the permutation $\sigma =\pi_1 \cdots \pi_{i-1}\boxed{\pi_j \pi_{j+1}\cdots \pi_{k-1}}\boxed{\pi_i \pi_{i+1}\cdots \pi_{j-1}} \pi_k \cdots \pi_n$.
\item As for the reversal, also for the transposition model there is a ``prefix variant".
In the \emph{prefix transposition model} the leftmost block of elements to interchange is a prefix of the permutation.
Sorting by prefix transposition is studied in \cite{DM}.
\end{itemize}

Independently from the chosen model, there are some general questions that can be asked in order to gain a better understanding of its combinatorial properties.
First of all, the operations of the model often (but not always) allow to define a \emph{distance} $d$ between two permutations $\rho$ and $\sigma$,
as the minimum number of elementary operations needed to transform $\rho$ into $\sigma$.
Moreover, when the operations are nice enough, the above distance $d$ could even be \emph{left-invariant},
meaning that, given permutations $\pi ,\rho ,\sigma$ (of the same length), $d(\pi ,\rho )=d(\sigma \pi ,\sigma \rho )$.
As a consequence, choosing for instance $\sigma =\rho^{-1}$,
the problem of evaluating the distance $d(\pi ,\rho )$ reduces to that of sorting $\pi$ with the minimum number of elementary allowed operations.
Now, if $d$ is a left-invariant distance on the set $S_n$ of all permutations of the same length,
define the \emph{$k$-ball} of $S_n$ to be the set $B_k ^{(d)}(n)=\{ \rho \in S_n \; |\; d(\rho ,id_n )\leq k\}$,
where $id_n$ is the identity permutation of length $n$.
The following questions are quite natural to ask:
\begin{itemize}
\item compute the diameter of $B_k ^{(d)}(n)$, i.e. the maximum distance between two permutations of $B_k ^{(d)}(n)$;
\item compute the diameter of $S_n$, i.e. the maximum distance between two permutations of $S_n$;
\item characterize the permutations of $\partial B_k ^{(d)}(n)$, i.e. the permutations of $B_k ^{(d)}(n)$ having maximum distance from the identity;
\item characterize the permutations of $\partial S_n$, i.e. the permutations of $S_n$ having maximum distance from the identity;
\item characterize and enumerate the permutations of $B_k ^{(d)}(n)$;
\item design sorting algorithms and study the related complexity issues.
\end{itemize}
In the literature there are several results, concerning several evolution models, which give some insight into the above problems.
Our work starts from the observation that, in many cases, the balls $B_k ^{(d)}(n)$ can be characterized in terms of \emph{pattern avoidance}.
Recall that, given two permutations $\sigma \in S_k$ and $\tau =\tau_1 \tau_2 \cdots \tau_n \in S_n$, with $k\leq n$,
we say that $\sigma$ is a \emph{pattern} of $\tau$ when
there exist $1\leq i_1 <i_2 <\cdots <i_k \leq n$ such that $\tau_{i_1}\tau_{i_2}\cdots \tau_{i_k}$ (as a permutation) is isomorphic to $\sigma$
(which means that $\tau_{i_1}, \tau_{i_2}, \ldots , \tau_{i_k}$ are in the same relative order as the elements of $\sigma$).
This notion of pattern in permutations defines a partial order, and the resulting poset is known as the \emph{permutation pattern poset}.
When $\sigma$ is not a pattern of $\tau$, we say that $\tau$ \emph{avoids} $\sigma$.
A down-set $I$ (also called a permutation class) of the permutation pattern poset can be described in terms of its minimal excluded permutations
(or, equivalently, the minimal elements of the complementary up-set): these permutations are called the \emph{basis} of $I$.
The idea of studying the balls $B_k ^{(d)}(n)$ in terms of pattern avoidance is not new.
As far as we know, the first model which has been investigated from this point of view is the (whole) tandem duplication-random loss model:
Bouvel and Rossin \cite{BR} have in fact shown that, in such a model,
the ball $B_k ^{(d)}=\bigcup_{n\geq 0}B_k ^{(d)}(n)$ is a class of pattern avoiding permutations,
whose basis is the set of minimal permutations having $d$ descents (here minimal is intended in the permutation pattern order).
Subsequent works \cite{BoPe,BF,CGM} have been done concerning the enumeration of the basis permutations of such classes.
More recently, Homberger and Vatter \cite{HV} described an algorithm for the enumeration of any polynomial permutation class,
which can be fruitfully used for all the above mentioned distances, since the resulting balls are indeed polynomial classes.
However, their results do not allow to find information on the basis of the classes.

In the present work we try to enhance what have been obtained in \cite{HV} in two directions.
First, we aim at giving a structural characterization of the balls for some of the above distances,
thus complementing the results in \cite{HV}, which is more concerned with computational issues.
Secondly, we provide some insight on the properties of the bases of such balls,
hoping to gain a better understanding of them.
In particular, we will prove that the problem of determining what we call the \emph{clean compact peg basis} of the balls
(in both the reversal and the prefix reversal models) is decidable.
The knowledge of the clean compact peg basis is of great help in finding the standard basis, as we will show in the next pages.

The companion paper \cite{CF}, published in the proceedings of the conference GASCom 2018,
deals with the cases of the block transposition and prefix block transposition models.
In the present work we find analogous results for the reversal and prefix reversal model.

\section{The reversal model}

\subsection{Generating permutations}

We start with the study of the general reversal model.
Recall (from the Introduction) that the reversal model is defined as the set of all permutations (of finite length)
endowed with a combinatorial operation, called \emph{reversal},
which consists of reversing a given segment (i.e. a set of adjacent elements in one-line notation) of a permutation,
thus obtaining a new permutation of the same length.
This operation defines a distance, called \emph{reversal distance},
which is defined as the minimum number of reversals needed to transform one permutation into another.
Since the reversal distance is left-invariant,
we can limit ourselves to computing the distance of a given permutation from the identity permutation
(having the correct length).
Formally, such a quantity will be denoted $rd(\pi)$ (for a given permutation $\pi$),
and the set of all permutations having distance $k$ from the identity (that is, the \emph{ball of radius $k$})
will be denoted $B_k^{(rd)}$.

\bigskip

In order to deal with the reversal model,
it will be useful to work both with standard permutations and with permutations whose elements are suitably decorated.
Following \cite{HV}, a \emph{peg permutation} of length $n$ is a decorated permutation
$\pi ^{\varepsilon}=\pi_1 ^{\varepsilon_1 }\pi_2 ^{\varepsilon_2 }\cdots \pi_n ^{\varepsilon_n }$,
where each $\varepsilon_i$ is either $+,-$ or $\bullet$; for instance, $3^+ 4^{\bullet}1^- 5^-2^+$ is a peg permutation of length 5.
Given a peg permutation $\pi^{\varepsilon}=\pi_1 ^{\varepsilon_1 }\pi_2 ^{\varepsilon_2 }\cdots \pi_n ^{\varepsilon_n }$,
we denote with $\pi$ the underlying permutation $\pi=\pi_1 \cdots \pi_n$.
There is a particular subclass of peg permutations that will be useful in what follows.
To define it, we need to introduce a few terminology.

An \emph{increasing strip} of $\pi ^{\varepsilon}$ is a maximal consecutive substring $\pi_i ^{\varepsilon_i } \cdots \pi_{i+k-1}^{\varepsilon_{i+k-1}}$
such that, for all $j$, $\pi_{j+1}=\pi_j +1$ and $\varepsilon_j$ is either $+$ or $\bullet$.
Similarly, a \emph{decreasing strip} of $\pi ^{\varepsilon}$ is a maximal consecutive substring
$\pi_i ^{\varepsilon_i } \cdots \pi_{i+k-1}^{\varepsilon_{i+k-1}}$
such that, for all $j$, $\pi_{j+1}=\pi_j -1$ and $\varepsilon_j$ is either $-$ or $\bullet$.
When a sequence of consecutive entries is either an increasing or a decreasing strip, we will simply call it a \emph{strip}.
An \textit{identity peg permutation} of length $n$ is an identity permutation in which every element is decorated with either $+$ or $\bullet$;
notice that there are $2^n$ such permutations.
Following again the terminology of \cite{HV}, a \emph{clean compact peg permutation} is a peg permutation all of whose strips have length 1.
A weaker notion is that of a \emph{compact peg permutation},
which is a peg permutation all of whose strips either have length 1 or are formed exclusively by elements decorated with $\bullet$.
Thus, for instance, the above peg permutation $3^+ 4^{\bullet} 1^- 5^- 2^+$ is not clean compact,
whereas $3^{\bullet} 4^{\bullet} 1^- 5^- 2^+$ is compact, and $2^+ 5^- 4^+ 1^{\bullet} 3^-$ is a clean compact peg permutation of length 5.
Any permutation $\pi$ can be associated with a unique clean compact peg permutation, denoted $peg(\pi)$,
which is obtained by replacing each strip of $\pi$ with its minimum element decorated in a specific way, then suitably rescaling the resulting word.
The decorations are the following:
\begin{enumerate}
\item if the element is the minimum of an increasing strip of length $\geq 2$, we decorate it with $+$;
\item if the element is the minimum of a decreasing strip of length $\geq 2$, we decorate it with $-$;
\item if the element is the unique element of a strip of length $1$, we decorate it with $\bullet$.
\end{enumerate}
For instance, if $\pi =32451678$, then $peg(\pi )=2^- 3^+ 1^{\bullet}4^+$.

\bigskip

We can extend the notion of reversal distance and pattern involvement to peg permutations.
First of all, if we want to apply a reversal to a peg permutation $\pi^{\varepsilon}$, we perform the standard reversal to the underlying permutation $\pi$,
then we invert the decoration $+$ and $-$ of each element involved in the operation (the $\bullet$ decorations remain unchanged).
A similar notion of \textit{oriented reversal} has already been considered in some classical works (see \cite{HP}).
For example, if we apply the reversal of indices $2,4$ to the peg permutation $3^+ 1^+ 2^- 5^{\bullet} 4^+$,
we obtain $3^+ \underbrace{1^+ 2^- 5^{\bullet}}_{rev} 4^+ \ = \ 3^+ 5^{\bullet} 2^+ 1^- 4^+$.
In analogy with the usual definition,
we will say that the reversal distance of a peg permutation $\pi^{\varepsilon}$ of length $n$ (from the identity)
is the minimum number of (oriented) reversals
needed to transform $\pi^{\varepsilon}$ in the corresponding identity permutation of length $n$.
From now on, we will use the same term ``reversal" for both standard and peg permutations,
omitting the word ``oriented" in the peg case.
Moreover, we will use the same symbol $rd$ to denote the reversal distance for of peg permutations.
The set of the peg permutations having distance at most $k$ from an identity permutation will be denoted with
$\hat{B}_k^{(rd)}$.

Concerning the notion of pattern containment, given two peg permutations
$\sigma^{\varepsilon}=\sigma_1^{\varepsilon_1} \cdots \sigma_k^{\varepsilon_k},\tau^{\varepsilon}=\tau_1^{\varepsilon_1} \cdots \tau_n^{\varepsilon_n}$,
with $k \le n$, we will say that $\sigma^{\varepsilon}$ is a \emph{pattern} of $\tau^{\varepsilon}$ when
there exist $1\leq i_1 <i_2 <\cdots <i_k \leq n$ such that:
\begin{enumerate}
\item $\tau_{i_1}\tau_{i_2}\cdots \tau_{i_k}$ (as a permutation) is isomorphic to $\sigma$;
\item if $\tau_{i_j}$ is decorated $+$ (resp. $-$), then $\sigma_j$ is decorated $+$ (resp. $-$).
\end{enumerate}

Note that, for example, $1^+ 2^{\bullet} 3^+$ is a pattern of $1^+ 2^- 3^+$, while the converse is not true.
The above relation is a partial order on the set of all peg permutations,
and the resulting poset is called the \emph{peg permutation pattern poset}.
We have the following fundamental result, whose easy proof is left to the reader.
\begin{prop}\label{Downset}
\begin{enumerate}
\item If $\sigma$ is a pattern of $\tau$ and $\tau \in B_k^{(rd)}$, then $\sigma \in B_k^{(rd)}$;
in other words, $B_k^{(rd)}$ is a down-set of the permutation pattern poset.
\item If $\sigma^{\varepsilon_s}$ is a pattern of $\tau^{\varepsilon_t}$
and $\tau^{\varepsilon_t} \in \hat{B}_k^{(rd)}$,
then $\sigma^{\varepsilon_s} \in \hat{B}_k^{(rd)}$;
in other words, $\hat{B}_k^{(rd)}$ is a down-set of the peg permutation pattern poset.
\end{enumerate}
\end{prop}

Now, given a peg permutation $\pi ^{\varepsilon}=\pi_1 ^{\varepsilon_1 }\pi_2 ^{\varepsilon_2 }\cdots \pi_n ^{\varepsilon_n }$,
we can create a set of associated permutations by monotone inflations.
A vector $v=(v_1 ,\ldots ,v_n)$ of nonnegative integers is a \textit{legal inflation vector} for $\pi^{\varepsilon}$ whenever
$v_i \in \{ 0,1 \}$ for each $i$ such that $\varepsilon_i=\bullet$;
then, for a peg permutation $\pi^{\varepsilon}$ and a legal inflation vector $v$,
we define the \textit{monotone inflation} of $\pi^{\varepsilon}$ through $v$ as
the permutation $\pi^{\varepsilon} [v]$ obtained from $\pi^{\varepsilon}$ in the following way:
\begin{enumerate}
\item if an element $\pi_i^{\varepsilon_i}$ is decorated $+$, we replace it with the identity permutation of length $v_i$;
\item if an element $\pi_j^{\varepsilon_j}$ is decorated $-$, we replace it with the reverse identity permutation of length $v_j$;
\item if an element $\pi_k^{\varepsilon_k}$ is decorated $\bullet$, we delete the element if $v_k=0$, otherwise we simply delete the sign $\varepsilon_k$;
\item finally, we suitably rescale each new insertion so to mantain the relative order of the elements of the permutation underlying $\pi^{\varepsilon}$.
\end{enumerate}

If we perform the same construction, but we inflate using identity and reverse identity peg permutations,
we obtain the \textit{peg monotone inflation} of $\pi^{\varepsilon}$ through $v$, denoted $\pi^{\varepsilon} [v]_{peg}$;
note that, in this case, we obtain a set of peg permutations for each choice of $v$.
For instance, if $\pi^{\varepsilon}=3^+ 1^+ 2^{\bullet} 5^- 4^{\bullet}$ and $v=(2,0,1,3,1)$
(note that $v$ is a legal inflation vector for $\pi^{\varepsilon}$), then we have that:
\begin{enumerate}
\item $\pi^{\varepsilon} [v] = \underbrace{23}_{3^+} \underbrace{\ldots}_{1^+} \underbrace{1}_{2^{\bullet}} \underbrace{765}_{5^-} \underbrace{4}_{4^{\bullet}}$.
\item $\pi^{\varepsilon} [v]_{peg} \ni \underbrace{2^{\bullet} 3^+}_{3^+} \underbrace{\ldots}_{1^+} \underbrace{1^{\bullet}}_{2^{\bullet}} \underbrace{7^- 6^{\bullet} 5^-}_{5^-} \underbrace{4^{\bullet}}_{4^{\bullet}}$.
\end{enumerate}

From now on, again following \cite{HV}, we will denote with $Grid(\pi^{\varepsilon})$ the \textit{grid class} of $\pi^{\varepsilon}$,
i.e. the set of all the monotone inflations of $\pi^{\varepsilon}$,
and with $Grid(C)$ the set $\bigcup_{\pi^{\varepsilon} \in C} Grid(\pi^{\varepsilon})$, for a given set $C$ of peg permutations.
The analogous classes for peg permutations will be denoted $Grid_{peg}(\pi^{\varepsilon})$ and $Grid_{peg}(C)$.
The notion of grid class is compatible with the pattern involvement relation, in the sense specified by the following proposition.
We give only a sketch of the proof, leaving the details to the reader.

\begin{prop}\label{grid_prop}
\begin{enumerate}
\item If $\pi^{\varepsilon}$ is a pattern of $\gamma^{\varepsilon'}$,
then $Grid(\pi^{\varepsilon}) \subseteq Grid(\gamma^{\varepsilon'})$
and $Grid_{peg}(\pi^{\varepsilon}) \subseteq Grid_{peg}(\gamma^{\varepsilon'})$.
\item If $\pi^{\varepsilon} \in \hat{B}_k^{(rd)}$,
then $Grid(\pi^{\varepsilon}) \subseteq B^{(rd)}_k$ and $Grid_{peg}(\pi^{\varepsilon}) \subseteq \hat{B}_k^{(rd)}$.
\end{enumerate}
\end{prop}

\emph{Proof.} \quad Given an occurrence of the pattern $\pi^{\varepsilon}$ in $\gamma^{\varepsilon'}$,
every monotone inflation of $\pi^{\varepsilon}$ can be obtained by
inflating the elements of $\gamma^{\varepsilon}$ that form that occurrence of $\pi^{\varepsilon}$ exactly in the same way and the remaining ones with $0$,
so the first statement follows.
To prove the second statement, observe that, given a sorting sequence of reversals for $\pi^{\varepsilon}$,
one can find a new sorting sequence of the same length for every permutation $\gamma$ in $Grid(\pi^{\varepsilon})$:
for each reversal $R$ of the given sorting sequence,
define a reversal $R'$ for $\gamma$ by incorporating in it all elements of $\gamma$ obtained by inflating each element of $\pi^{\varepsilon}$ affected by $R$.
An analogous argument works also in the case of grid peg classes.\cvd

\begin{cor}\label{peg_dist}
Let $\pi=\pi_1 \cdots \pi_n$
and $peg(\pi)=\alpha_1^{a_1} \cdots \alpha_k^{a_k}$ the associated clean compact peg permutation.
Then $rd(\pi) \le rd(peg(\pi))$.
\end{cor}

\emph{Proof.} \quad By construction, we have that $\pi \in Grid(peg(\pi))$,
so $rd(\pi) \le rd(peg(\pi))$ from the second statement of the above proposition.
\cvd

Note that the opposite inequality is not true in general.
For example, $rd(3412)=2$
(the sequence of reversals is $\underbrace{341}2 \rightsquigarrow 1\underbrace{432} \rightsquigarrow 1234$),
however the associated clean compact peg permutation is $2^+ 1^+$, that cannot be sorted with $2$ reversals.

The notions of monotone inflation and peg monotone inflation allow us to give
a nice and precise description of $B_k ^{(rd)}=\bigcup_{n\geq 0} B_k ^{(rd)}(n)$.
Our goal is to prove that $B_k^{(rd)}$ is
the union of the grid classes of the clean compact permutations of $B_k^{(rd)}$
that are maximal with respect to the peg pattern involvement relation.
To do this, we need some additional results and constructions.

\begin{lemma}\label{max_no_bullet}
Let $\pi^{\varepsilon}$ be a clean compact peg permutation in $\hat{B}_k^{(rd)}$
that is maximal with respect to the peg pattern involvement relation;
then $\pi^{\varepsilon}$ has no $\bullet$ decoration.
\end{lemma}

\emph{Proof.} \quad Suppose there is an index $j$ such that $\pi_j$ is decorated $\bullet$ in $\pi^{\varepsilon}$.
Given a sequence of (at most) $k$ reversals that sorts $\pi^{\varepsilon}$,
let $t$ be the number of reversals that involve the element $\pi^{\bullet}_j$.
If $t$ is even, then the same sequence of reversals also sorts the permutation $\bar{\pi}^{\varepsilon}$,
obtained by decorating $\pi_j$ with $+$ instead of $\bullet$.
This means that $\bar{\pi}^{\varepsilon}$ has distance (at most) $k$ and contains $\pi^{\varepsilon}$ as a pattern;
moreover, $\bar{\pi}^{\varepsilon}$ is clean compact,
because $\pi^{\varepsilon}$ is clean compact and we are replacing a $\bullet$ with a $+$.
Thus we obtain a contradiction since $\pi^{\varepsilon}$ is supposed to be maximal.
If $t$ is odd, we can argue in the same way by replacing the $\bullet$ with a $-$ decoration.
\cvd

Suppose $\pi^{\varepsilon}=\pi_1^{\varepsilon_1} \cdots \pi_n^{\varepsilon_n}$ is a peg permutation of length $n$,
with no elements decorated $\bullet$.
Inflate $\pi^{\varepsilon}$ by choosing two (not necessarily distinct) indices $1 \leq i \leq j \leq n$
and replacing $\pi_i^{\varepsilon_i}$ and $\pi_j^{\varepsilon_j}$ as follows:
\begin{enumerate}
\item if $i \neq j$, then replace $\pi_i$ with $\pi_i ^+ (\pi_i+1)^+$, if $\varepsilon_i=+$,
or with $(\pi_i+1)^- \pi_i^-$ if $\varepsilon_i=-$; do the same for $\pi_j$;
\item if $i=j$, then replace $\pi_i$ with $\pi_i^+ (\pi_i+1)^+ (\pi_i+2)^+$, if $\varepsilon_i=+$;
otherwise, if $\varepsilon_i=-$, replace it with $(\pi_i+2)^- (\pi_i+1)^- \pi_i^-$.
\end{enumerate}

Now rescale the resulting string according to the relative order of the elements of $\pi$,
so to obtain a new peg permutation of length $n+2$.
If $I=\{i,j \}$, $i \le j$, is the multiset of the selected indices,
the resulting peg permutation will be denoted $\pi^{\varepsilon}_I$.
A simple case by case analysis shows that, if the starting permutation $\pi^{\varepsilon}$ is clean compact,
there exists a unique reversal $\rho_I$ for $\pi^{\varepsilon}_I$ such that the resulting permutation
is clean compact as well;
more specifically, $\rho_I$ is the reversal with indices $i+1,j+1$.
Call $\tilde{\pi}^{\varepsilon}_I$ the permutation obtained by applying $\rho_I$ to $\pi^{\varepsilon}_I$.
As an example, consider the peg permutation $\pi^{\varepsilon} =1^+ 2^- 3^+$ and the multiset of indices $I=\{ 1,3 \}$;
then we get $\pi^{\varepsilon}_I =1^+ 2^+ \ 3^- \ 4^+ 5^+$ and
$\tilde{\pi}^{\varepsilon}_I =1^+ \mathbf{4^- 3^+ 2^-} 5^+$.
Choosing instead $J=\{ 2,2 \}$, we get $\pi^{\varepsilon}_J=1^+ \ 4^- 3^- 2^- \ 5^+$
and $\tilde{\pi}^{\varepsilon}_J=1^+ 4^- \mathbf{3^+} 2^- 5^+$.
Note that in the two above cases we get the same resulting permutation;
this shows how, in general, different choices of multisets can lead to the same permutation.

The following lemma sums up a few basic properties of the above construction that will be useful in the sequel.
We leave the easy proof to the reader.

\begin{lemma}\label{prel}
Let $\pi^{\varepsilon} =\pi_1^{\varepsilon_1} \cdots \pi_n^{\varepsilon_n}$ be
a clean compact peg permutation of length $n$ without $\bullet$,
and $I$ a multiset of indices of $\pi$ of cardinality $2$.
Then $\tilde{\pi}^{\varepsilon}_I$ is a clean compact peg permutation of length $n+2$ without $\bullet$.
Moreover, if $\pi^{\varepsilon}_1 =1^+$ and $\pi^{\varepsilon}_n =n^+$,
then $(\tilde{\pi}^{\varepsilon}_I )_1 =1^+$ and $(\tilde{\pi}^{\varepsilon}_I )_{n+2} =(n+2)^+$.
\end{lemma}

We are now ready to state the proposition that provides an inductive description of $B_k ^{(rd)}$.
The (excessively technical) proof is omitted.

\begin{prop}\label{plusone}
Let $\mathcal{I}(n)$ be the set of all multisets of cardinality $2$ of $\{ 1,2,\ldots n\}$
and let $\pi^{\varepsilon}$ be a clean compact peg permutation of length $n$ without $\bullet$.
Denote with $Grid(\pi^{\varepsilon})^{+1}$ the set of all permutations that can be obtained
by applying a single reversal to each permutation of $Grid(\pi^{\varepsilon})$
and with $Grid_{peg}(\pi^{\varepsilon})^{+1}$ the corresponding set of peg permutations.
Then:
\begin{enumerate}
\item $Grid(\pi^{\varepsilon})^{+1}=\bigcup_{I \in \mathcal{I}(n)}Grid(\tilde{\pi}^{\varepsilon}_I )$;
\item $Grid_{peg}(\pi^{\varepsilon})^{+1}=\bigcup_{I \in \mathcal{I}(n)}Grid_{peg}(\tilde{\pi}^{\varepsilon}_I )$.
\end{enumerate}
\end{prop}

The above proposition tells that
each permutation obtained by applying a single reversal to a permutation in $Grid(\pi^{\varepsilon})$
is a monotone inflation of $\tilde{\pi}^{\varepsilon}_I$, for some multiset $I$.
An analogous fact holds in the poset of peg permutations.

A repeated applications of the above proposition leads to the following result.

\begin{theorem}\label{Gener}
For every $k \ge 1$,
there exist $N=N(k)$ clean compact peg permutations
$\alpha^{\varepsilon^1}_{(1)},\ldots ,\alpha^{\varepsilon^N}_{(N)}$,
of length $2k+1$ and having distance $k$ from the identity, such that:
\begin{enumerate}
\item $B_k ^{(rd)}=\bigcup_{j=1}^{N}Grid(\alpha^{\varepsilon^j}_{(j)})$;
\item $\hat{B}_k^{(rd)}=\bigcup_{j=1}^{N}Grid_{peg}(\alpha^{\varepsilon^j}_{(j)})$.
\end{enumerate}
Such permutations will be called \emph{k-generating permutations}
and the set $\{ \alpha_{(j)}^{\varepsilon_j}\, |\, j=1,\ldots ,N\}$
is the $k$-\emph{generating set} of $B_k^{(rd)}$ and $\hat{B}_k^{(rd)}$.
\end{theorem}

\emph{Proof.}\quad
\begin{enumerate}
\item We proceed by induction on $k$. Obviously $B_0^{(rd)}=Grid(1^+)$.
When $k=1$, it is easy to observe that $B_1 ^{(rd)}=Grid(1^+ 2^- 3^+)$,
and $1^+ 2^- 3^+$ is a clean compact peg permutation, has length $3$ and has distance $1$ from the identity.
Now consider a permutation $\overline{\pi}\in B_{k+1} ^{(rd)}\setminus B_k ^{(rd)}$;
this means, in particular, that there is a permutation $\pi \in B_k ^{(rd)}$ such that
$\overline{\pi}$ is obtained from $\pi$ by applying a single reversal.
Thus, using the induction hypothesis,
we can assert that there exists a clean compact peg permutation $\alpha^{\varepsilon}$,
having length $2k+1$ and distance $k$ from the identity, such that $\pi \in Grid(\alpha^{\varepsilon})$.
By Proposition \ref{plusone},
there exists $I \in \mathcal{I}(2k+1)$ such that $\overline{\pi} \in Grid(\tilde{\alpha}^{\varepsilon}_I )$.
Notice that, by Lemma \ref{prel}, $\tilde{\alpha}^{\varepsilon}_I$ is a clean compact peg permutation of length $2k+3$.

It remains to prove that $\tilde{\alpha}^{\varepsilon}_I$ has distance $k+1$ from the identity.
Clearly $rd(\tilde{\alpha}^{\varepsilon}_I )\leq k+1$.
Conversely,
define a \textit{breakpoint} of a peg permutation to be a pair of adjacent elements that is not part of a strip.
In particular, analogously to what happens for standard permutations \cite{HP},
a single reversal can change the number of breakpoints by at most $2$.
Since $\tilde{\alpha}^{\varepsilon}_I$ is clean compact and has length $2k+3$, it has exactly $2k+2$ breakpoints,
whereas an identity peg permutation of length $2k+3$ has no breakpoint.
Therefore at least $\frac{2k+2}{2}=k+1$ reversals are needed to sort $\tilde{\alpha}^{\varepsilon}_I$.
This means that $rd(\tilde{\alpha}^{\varepsilon}_I) \ge k+1$, as desired.

\item The proof is identical to that of the previous case.
Just notice that $\hat{B}_0^{(rd)}= \left\lbrace 1^+, 1^{\bullet} \right\rbrace = Grid_{peg}(1^+)$,
so the base case of the induction argument is the same as before.\hfill $\blacksquare$
\end{enumerate}

Notice that the maximum length of a clean compact peg permutation of $\hat{B}_k^{(rd)}$ is $2k+1$,
since each reversal can create at most $2$ new breakpoints and an identity permutation has no breakpoint.
Therefore, recalling that the permutations $\alpha^{\varepsilon^1}_{(1)},\ldots ,\alpha^{\varepsilon^N}_{(N)}$ do not have elements decorated $\bullet$,
we have that each generating permutation is maximal inside $\hat{B}_k^{(rd)}$ (because it has maximum length and no decoration $\bullet$)
and clean compact (by construction).
The next theorem proves the converse.

\begin{theorem}
Let $k \ge 1$ and let $\pi^{\varepsilon}$ be a maximal clean compact peg permutation in $\hat{B}_k^{(rd)}$.
Then $\pi^{\varepsilon}$ is a $k$-generating permutation.
\end{theorem}

\emph{Proof.} \quad The maximality of $\pi^{\varepsilon}$ guarantees that it has no element decorated $\bullet$;
also notice that $\pi^{\varepsilon}$ has length $\le 2k+1$, because it is clean compact and belongs to $\hat{B}_k^{(rd)}$.
Moreover, Theorem \ref{Gener} implies that there exists a generating permutation $\alpha^{\varepsilon_a}$,
of length $2k+1$ and having distance $k$ from the identity, such that $\pi^{\varepsilon} \in Grid_{peg}(\alpha^{\varepsilon_a})$,
i.e. there is a vector $v$ of nonnegative integers such that $\pi^{\varepsilon} \in \alpha^{\varepsilon_a} [v]_{peg}$.
In particular, it must be $v_i \le 1$ for each $i$, again because $\pi^{\varepsilon}$ is clean compact.
As a consequence, we have that $\pi^{\varepsilon}$ is a pattern of $\alpha^{\varepsilon_a}$,
but $\pi^{\varepsilon}$ is a maximal clean compact element of $\hat{B}_k^{(rd)}$, so it has to be $\pi^{\varepsilon}=\alpha^{\varepsilon_a}$, as desired.
\cvd

\begin{cor}\label{genblock}
For every $k \ge 1$, $B_k^{(rd)}$ is the union of the grid classes of the maximal clean compact peg permutations of $\hat{B}_k^{(rd)}$.
Thus $\alpha^{\varepsilon}$ is a $k$-generating permutations if and only if:
\begin{enumerate}
\item $rd(\alpha^{\varepsilon})=k$;
\item $\alpha^{\varepsilon}$ has length $2k+1$;
\item $\alpha^{\varepsilon}$ is clean compact;
\item $\alpha^{\varepsilon}$ has no decorations $\bullet$.
\end{enumerate}
\end{cor}

The previous results suggest a procedure to list the generating permutations of $B_k ^{(rd)}$:
starting from $1^+$,
perform $k$ successive monotone inflations as in Lemma \ref{prel} so to obtain all generating permutations of $B_k ^{(rd)}$.
This is similar to the approach used in \cite{HV}.
For instance, when $k=2$, the generating set for $B_2 ^{(rd)}$ consists of the four permutations
$1^+ 2^- 3^+ 4^- 5^+, 1^+ 4^- 3^+ 2^- 5^+, 1^+ 4^+ 2^- 3^- 5^+, 1^+ 3^- 4^- 2^+ 5^+$.
As we have already observed, it is possible to obtain the same generating permutation several times.
This is the main reason for which this approach cannot be used to enumerate the generating set.

\bigskip

\textbf{Open problem 1.}\quad \emph{Enumerate the generating permutations of $B_k ^{(rd)}$, for every $k$}
(that is, determine the quantity $N(k)$ in the statement of Theorem \ref{Gener}).

\subsection{The bases}

As we have observed in Proposition \ref{Downset}, $B_k^{(rd)}$ and $\hat{B}_k^{(rd)}$ are down-sets
in the permutation pattern poset and in the peg permutation pattern poset, respectively.
Therefore it would be nice to find some information concerning their bases.

We consider the ball $\hat{B}_k^{(rd)}$ first and try to find some properties of its basis.
For technical reasons that will become evident soon,
in the sequel we restrict our attention to the poset of clean compact permutations as a subposet of the whole peg permutation pattern poset.
Here we investigate properties of the set of minimal excluded permutations of $\hat{B}_k^{(rd)}$,
which can be called the \emph{clean compact peg basis} of $\hat{B}_k^{(rd)}$.
Our goal is to determine an upper bound for the length of a basis permutation;
such a bound can be used as a stopping condition for an algorithm which tries to explicitly describe the clean compact peg basis of $\hat{B}_k^{(rd)}$.
This shows, in particular, that the problem of finding the clean compact peg basis is decidable.
We start by proving some useful facts.

\begin{lemma}\label{Reduced_pattern}
Let $\pi^{\varepsilon}=\pi_1^{\varepsilon_1} \cdots \pi_n^{\varepsilon_n}$ be a clean compact peg permutation;
then $\pi^{\varepsilon}$ contains at least one pattern $\gamma^{\varepsilon'}$ of length $n-1$ which is clean compact as well.
\end{lemma}

\emph{Proof.} \quad Consider the index $i$ such that $\pi_i=n$.
If $i=1$ or $i=n$, then we can remove $\pi_i^{\varepsilon_i}$, thus obtaining a peg permutation that is still clean compact.
Otherwise we need to consider a few distinct cases.
Below we give details only for one of them, leaving the (analogous) remaining proofs to the reader.
More specifically, we provide a complete analysis of the case in which $n$ is decorated $+$,
and skip the cases in which $n$ is decorated $-$ or $\bullet$.

Suppose that $\pi_i^{\varepsilon_i}=n^+$.
In most of the situations we can remove $n^+$ and obtain a clean compact peg permutations.
There are only two ``bad"  cases,
occurring when the two elements on the left and on the right of $n^+$ (respectively) form
an increasing strip or a decreasing strip (of length 2).
We now analyze only the first of the two above cases, the other one being completely analogous.
Our hypothesis is therefore that $\pi^{\varepsilon}$ contains a string (of adjacent elements) of the form
$k^{+ / \bullet} \ n^+ \ (k+1)^{+ / \bullet}$, for some $k$.
Here the decoration $+ / \bullet$ stands for either a $+$ or a $\bullet$.
In such a situation, removing $k+1$ almost always leads to a clean compact peg permutation,
apart from two bad cases.
\begin{enumerate}
\item Suppose that $\pi^{\varepsilon}$ contains a string of the form
$(k+2)^{- / \bullet} \ k^{\bullet} \ n^+ \ (k+1)^{+ / \bullet}$.
Then we can remove $k$, unless we are in one of the two cases described below,
which we will now manage separately.
\begin{itemize}
\item[(i)] The first bad case occurs when $n=k+3$, and so $\pi^{\varepsilon}$ contains a string of the form
$(k+2)^{\bullet} \ k^{\bullet} \ (k+3)^+ \ (k+1)^{+ / \bullet}$.
We can then remove $k+2$, which always works except when there is a string of the form
$(k-1)^{+ / \bullet} \ (k+2)^{\bullet} \ k^{\bullet} \ (k+3)^+ \ (k+1)^{+ / \bullet}$.
In such a case, however, we can always remove $k-1$ to obtain a clean compact peg permutation of length $n-1$.
\item[(ii)] The second bad case occurs when $\pi^{\varepsilon}$ contains a string of the form
$(k+2)^{- / \bullet} \ k^{\bullet} \ n^+ \ (k+1)^{\bullet} \ (k-1)^{- / \bullet}$.
But now we can remove $k-1$ and we are done.
\end{itemize}
\item Suppose that $n=k+2$, and so $\pi^{\varepsilon}$ contains a string of the form
$k^{+ / \bullet} \ (k+2)^+ \ (k+1)^{+ / \bullet}$.
So a good choice is to remove $k$, unless there is a string of the form
$k^{+ / \bullet} \ (k+2)^+ \ (k+1)^{\bullet} \ (k-1)^{- / \bullet}$.
In this case, however, one could remove $k-1$, unless there is a string of the form
$(k-2)^{+ / \bullet} \ k^{+ / \bullet} \ (k+2)^+ \ (k+1)^{\bullet} \ (k-1)^{- / \bullet}$.
Now it turns out that this argument can be repeated, so that we can always find an element to remove,
unless we have already considered all the elements of the permutation.
If this happens, then $\pi^{\varepsilon}$ has $1$ either as its first or its last element,
so $1$ can be removed.\cvd
\end{enumerate}

\begin{lemma}\label{basis_prop}
Let $\pi^{\varepsilon}=\pi_1^{\varepsilon_1} \cdots \pi_n^{\varepsilon_n}$ be a permutation
in the clean compact peg basis of $\hat{B}_{k}^{(rd)}$.
Then $\pi_1^{\varepsilon_1} \neq 1^{+ / \bullet}$ and $\pi_n^{\varepsilon_n} \neq n^{+ / \bullet}$.

\end{lemma}
\emph{Proof.}
If $\pi_1^{\varepsilon_1}=1^+$ or $\pi_1=1^{\bullet}$, then we could remove it without changing the reversal distance of $\pi^{\varepsilon}$ , against the minimality of $\pi^{\varepsilon}$. The case $\pi_n ^{\varepsilon_n}=n^{+ / \bullet}$ is analogous.
\cvd

We can now state our main result concerning the length of a basis permutation.

\begin{theorem}\label{basis}
Every permutation belonging to the clean compact peg basis of $\hat{B}_{k}^{(rd)}$ has length at most $2k+1$.
\end{theorem}

\emph{Proof.} \quad It is easy to prove that every basis permutation $\pi^{\varepsilon}$ has length at most $2k+2$.
Indeed, suppose that $\pi^{\varepsilon}$ has length (at least) $2k+3$.
By Lemma \ref{Reduced_pattern},
$\pi^{\varepsilon}$ contains a clean compact pattern $\gamma^{\varepsilon'}$ of length (at least) $2k+2$,
that must belong to $\hat{B}_{k}^{(rd)}$ (thanks to the minimality property of $\pi^{\varepsilon}$).
But the maximum length of a clean compact permutation in $\hat{B}_{k}^{(rd)}$ is $2k+1$, so we have a contradiction.

Now suppose that $\pi^{\varepsilon}$ has length $2k+2$;
again, $\pi^{\varepsilon}$ contains a clean compact pattern $\gamma^{\varepsilon'}$ of length $2k+1$, with $\gamma^{\varepsilon'} \in \hat{B}_{k}^{(rd)}$.

By Theorem \ref{Gener}, 
there exists a $k$-generating permutation $\alpha^a$ such that $\gamma^{\varepsilon'} \in \alpha^a [v]$,
for some legal inflation vector $v$.
In particular, $v_i \le 1$ for each $i$, because $\gamma^{\varepsilon'}$ is clean compact;
moreover, both $\gamma^{\varepsilon'}$ and $\alpha^a$ have length $2k+1$, so it must be $v_i =1$, for each $i$.
As a consequence of Lemma \ref{prel}, 
we have that $\gamma_1^{\varepsilon'_1}=1^{+ / \bullet}$ and $\gamma_{2k+1}^{\varepsilon'_{2k+1}}=(2k+1)^{+ / \bullet}$,
whereas Lemma \ref{basis_prop} implies that 
$\pi_1^{\varepsilon_1} \neq 1^{+ / \bullet}$ and $\pi_{2k+2} \neq (2k+2)^{+ / \bullet}$.
Since $\gamma^{\varepsilon'}$ is obtained from $\pi^{\varepsilon}$ by removing a single element, 
there are only two possible cases:
\begin{itemize}
\item[(i)] $\pi^{\varepsilon}=(2k+2)^{\varepsilon_1} 1^{+ / \bullet}\pi_3 ^{\varepsilon_3} \cdots \pi_{2k+1} ^{\varepsilon_{2k+1}}(2k+1)^{+ / \bullet}$
and $\gamma^{\varepsilon'}$ is obtained by removing the first element $2k+2$.
Then the permutation $\beta^{b}$ obtained from $\pi^{\varepsilon}$ by removing $1^{+ / \bullet}$ would be
a clean compact peg permutation of length $2k+1$ of $\hat{B}_{k}^{(rd)}$.
Therefore, again using Lemma \ref{prel}, it should be $\beta_1^{b_1}=1^{+ / \bullet}$, which gives a contradiction.
\item[(ii)] $\pi^{\varepsilon}=2^{+ / \bullet}\pi_2 ^{\varepsilon_2} \cdots \pi_{2k}^{\varepsilon_{2k}}(2k+2)^{+ / \bullet} 1^{\varepsilon_{2k+2}}$
 and $\gamma^{\varepsilon'}$ is obtained by removing the last element $1$.
 This case is clearly symmetric to the previous one: now we can remove $2k+2$ and use the same arguments.\cvd
\end{itemize}

\begin{cor}
The clean compact peg basis of $\hat{B}_{k}^{(rd)}$ is finite, for every $k \ge 1$.
\end{cor}

The above theorem also suggests a procedure to determine the clean compact peg basis of $\hat{B}_{k}^{(rd)}$
(thus showing that such a problem is decidable).
For each clean compact peg permutation $\pi^{\varepsilon}$ of length $2k+1$ which is not $k$-generating, take the set of permutations it covers:
%
\begin{itemize}
\item[-] if all of them lie below some $k$-generating permutation, then $\pi^{\varepsilon}$ is in the clean compact peg basis of $\hat{B}_{k}^{(rd)}$;
\item[-] otherwise, repeat the same procedure starting from the permutations covered by $\pi^{\varepsilon}$ which do not belong to $\hat{B}_{k}^{(rd)}$.
\end{itemize}

As an instance, we have the following result.

\begin{prop}
The clean compact peg basis of $\hat{B}_{1}^{(rd)}$ is $\{1^- 2^-, 2^+ 1^{\bullet}, 2^{\bullet} 1^+ \}$.
\end{prop}

\emph{Proof.}\quad Since $\hat{B}_{1}^{(rd)}=Grid_{peg}(1^+ 2^- 3^+)$,
we perform the above procedure with each clean compact peg permutations of length $3$, except for $1^+ 2^- 3^+$.
A direct inspection shows that there are no such permutations covering only elements of $\hat{B}_{1}^{(rd)}$.
Moreover $1^- 2^-, 2^+ 1^{\bullet}, 2^{\bullet} 1^+$ are the only clean compact peg permutations of length $2$ which are not in $\hat{B}_{1}^{(rd)}$
and such that all of their coverings are in $\hat{B}_{1}^{(rd)}$.\cvd

Unfortunately, translating the above results into the whole peg permutation pattern poset seems to be very hard.
The main problem is that there exist permutations in the peg basis which are not clean compact;
for example, the permutation $2^{\bullet} 1^{\bullet} 4^{\bullet} 3^{\bullet}$ has distance $2$ and is minimal with respect to this property:
it is in fact a peg basis element of $\hat{B}_{1}^{(rd)}$.
Also notice that both $2^{\bullet} 1^{\bullet} 4^{\bullet} 3^{\bullet}$ and $1^- 2^-$ are elements of the peg basis,
but $2^{\bullet} 1^{\bullet} 4^{\bullet} 3^{\bullet} \in Grid(1^- 2^-)$;
this suggests that there are some intrinsic difficulties in giving a satisfactory description of the peg basis,
which means that perhaps the whole peg permutation pattern poset is not the right setting for this problem.    
Moreover it seems to be quite hard to find a nontrivial upper bound for the length of the permutations of the peg basis,
although we know that such a basis is still finite as a consequence of some general results concerning grid classes (see \cite{AABRV}).
A useful characterization provided in \cite{HV} guarantees that a permutation of the peg basis is at least compact.

\begin{lemma}\cite{HV}
Let $\pi^{\varepsilon}$ be a peg permutation;
then $\pi^{\varepsilon}$ is compact if and only if $Grid(\pi^{\varepsilon})$ properly contains $Grid(\tau^t)$,
for each $\tau^t < \pi^{\varepsilon}$ in the peg permutation pattern poset.
\end{lemma}

\begin{prop}
Let $\pi^{\varepsilon}$ be a peg basis permutation of $\hat{B}_{k}^{(rd)}$. Then $\pi^{\varepsilon}$ is compact.
\end{prop}

\emph{Proof.} \quad Suppose that $\pi^{\varepsilon}$ contains a proper peg pattern $\tau^t$ such that $Grid(\tau^t) = Grid(\pi^{\varepsilon})$.
Then Proposition \ref{Downset} implies that $rd(\tau^t) \le rd(\pi^{\varepsilon})=k+1$.
Moreover it cannot be $rd(\tau^t) \le k$, otherwise we would also have $rd(\pi^{\varepsilon}) \le k$,
because $\pi^{\varepsilon} \in Grid(\tau^t)$ by hypothesis.
Thus $rd(\pi^{\varepsilon})=rd(\tau^t)=k+1$ and $\tau^t < \pi^{\varepsilon}$,
which gives a contradiction because $\pi^{\varepsilon}$ is minimal outside $\hat{B}_{k}^{(rd)}$.
Therefore, thanks to the above lemma, we can conclude that $\pi^{\varepsilon}$ is compact.\cvd

The final part of the present section is devoted to the standard basis of $B_k ^{(rd)}$.
The general theory of geometric grid classes allows to say that $B_k^{(rd)}$ is a permutation class having finite basis;
moreover $B_k^{(rd)}$ is \emph{strongly rational}, meaning that its generating function is rational,
together with the generating functions of all of its subclasses \cite{AABRV}.
Here we sketch a description the basis of $B_k^{(rd)}$, starting from the knowledge of the clean compact peg basis of $\hat{B}_{k}^{(rd)}$.

Let $k \ge 1$ and suppose that $\{ \beta_1^{b_1},\dots,\beta_N^{b_N} \}$ is the clean compact peg basis of $\hat{B}_{k}^{(rd)}$.
For every $i=1,\dots,N$, define the following sets of (standard) permutations:
\begin{itemize}
\item[-] $\mathcal{A}_{\beta_i^{b_i}} = \left\lbrace \pi: \ peg(\pi) = \beta_i^{b_i} \right\rbrace$;
\item[-] $\mathcal{M}_{\beta_i^{b_i}} = \left\lbrace \pi \in \mathcal{A}_{\beta_i^{b_i}}:
\pi \mbox{ is minimal such that } rd(\pi) = rd(\beta_i^{b_i}) \right\rbrace$;
\item[-] $\mathcal{M}=\bigcup_{i} \mathcal{M}_{\beta_i^{b_i}}$.
\end{itemize}

Our goal is to show that $B_k^{(rd)}=Av(\mathcal{M})$, so that the minimal elements in $\mathcal{M}$ form the basis of $B_k^{(rd)}$.
First of all, we prove that the sets $\mathcal{M}_{\beta_i^{b_i}}$ are not empty.

\begin{lemma}
Let $\pi^{\varepsilon}$ be a clean compact peg permutation such that $rd(\pi^{\varepsilon})=k$.
Then there is a permutation $\gamma \in Grid(\pi^{\varepsilon})$ such that
$rd(\gamma)=k$ and $peg(\gamma)=\pi^{\varepsilon}$.
\end{lemma}

\emph{Proof.} \quad Let $\gamma = \pi^{\varepsilon} [v]$, with $v$ defined as follows:
\begin{itemize}
\item[-] $v_i=1$, for each $i$ such that $\pi_i$ is decorated $\bullet$;
\item[-] $v_i=N>1$, otherwise.
\end{itemize}
Note that $v$ is a legal inflation vector for $\pi^{\varepsilon}$.
Suppose that $rd(\gamma) \le k-1$.
Theorem \ref{Gener} implies that there exists a $(k-1)$-generating permutation $\alpha^a$ such that
$\gamma \in Grid(\alpha^a)$.
In particular, $\alpha^a$ is a clean compact peg permutation of length $2k-1$.
Depending on the types of decorations $\pi^{\varepsilon}$ contains, we can distinguish two cases.
If $\pi^{\varepsilon}$ only contains $\bullet$,
the underlying permutation $\pi$ is $\gamma$ itself and $rd(\pi )=rd(\pi^{\varepsilon})=k$,
which is in contradiction with the assumption $rd(\gamma )\leq k-1$.
If instead $\pi^{\varepsilon}$ contains at least one decoration different from $\bullet$,
we observe that there is at least one element $\pi_j^{\varepsilon_j}$ of $\pi^{\varepsilon}$ such that
the corresponding strip $S_j$ in $\gamma$ is not contained in the inflation of a single element of $\alpha^a$:
otherwise, in fact, we could write $\pi^{\varepsilon}$ as an appropriate inflation of $\alpha^a$,
so it would be $\pi^{\varepsilon} \in Grid_{peg}(\alpha^a)$ and thus $rd(\pi) \le k-1$, which is not.
In particular, the decoration $\varepsilon_j$ is either $+$ or $-$,
because the corresponding strip $S_j$ must contain at least $2$ elements.
Moreover, if $S_j$ is increasing,
its elements are necessarily distributed among the inflations of entries of $\alpha^a$
which are decorated either $-$ or $\bullet$,
otherwise they could be assigned to the inflation of a single element decorated $+$
and so again $\pi^{\varepsilon} \in Grid_{peg}(\alpha^a)$.
An analogous argument holds if $S_j$ is decreasing.
Hence each element of $S_j$ must belong to a different strip of $\alpha^a$.
Now, recalling that $\alpha^a$ has length $|\alpha^a|=2k-1$, if we choose $N > 2k-1$,
we have that the $N$ elements of $S_j$
have to be distributed among the inflations of $N > |\alpha^a|$ different element of $\alpha^a$,
which is impossible.\cvd

In other words,
the above lemma says that the reversal distance of a peg permutation is equal to
the maximum distance of a permutation in its standard grid class:
every permutation in the grid class of $\pi^{\varepsilon}$ has distance at most $rd(\pi^{\varepsilon})$
(by Proposition \ref{grid_prop}) and we have proved that there is at least one permutation that attains this value.
As an example, consider the peg permutation $2^+ 1^+$, which has distance $3$.
The generating permutations of $\hat{B}_{2}^{(rd)}$ have length $5$, so we can choose $N=6>5$.
Define $\gamma=2^+1^+ \cdot [6,6]=\underbrace{789 \ 10 \ 11 \ 12}_{S_1} \ \underbrace{123456}_{S_2}$.
It cannot be $rd(\gamma)=2$,
since the six elements of either $S_1$ or $S_2$ should belong to
six distinct inflations of elements of a $2$-generating permutation, whose length is 5.


The next result is an immediate consequence of the above lemma.

\begin{cor}
If $\beta^b$ belongs to the clean compact peg basis of some $\hat{B}_{k}^{(rd)}$,
then the set $\mathcal{M}_{\beta^b}$ is not empty.
\end{cor}

We are now able to prove that $B_k^{(rd)}=Av(\mathcal{M})$.
First of all,
every permutation in $\mathcal{M}$ has distance equal to some permutation of the
clean compact peg basis of $\hat{B}_{k}^{(rd)}$;
therefore, if $\pi \in \mathcal{M}$, then $\pi \notin B_k^{(rd)}$.
Moreover,
we have to show that each permutation $\gamma \notin B_k^{(rd)}$ contains some permutation $\pi \in \mathcal{M}$
as a pattern.
If $\gamma \notin B_k^{(rd)}$, then also $peg(\gamma) \notin \hat{B}_{k}^{(rd)}$ by Corollary \ref{peg_dist}.
Since $peg(\gamma)$ is clean compact, there is a clean compact peg basis permutation $\beta_i^{b_i}$, for some $i$,
such that $peg(\gamma) \ge \beta_i^{b_i}$.
Given an occurrence of $\beta_i^{b_i}$ in $peg(\gamma)$,
consider the permutation $\hat{\gamma}$ obtained by taking only the strips of $\gamma$ corresponding
to the elements of $\beta_i^{b_i}$.
Clearly $\hat{\gamma} \le \gamma$ and moreover $peg(\hat{\gamma})=\beta_i^{b_i}$ by construction,
so $\hat{\gamma} \in \mathcal{A}_{\beta_i^{b_i}}$
and there exists a minimal permutation $\pi \in \mathcal{M}_{\beta_i^{b_i}}$ such that $\hat{\gamma} \ge \pi$.
Thus we have $\gamma \ge \hat{\gamma} \ge \pi$, as desired.

\bigskip

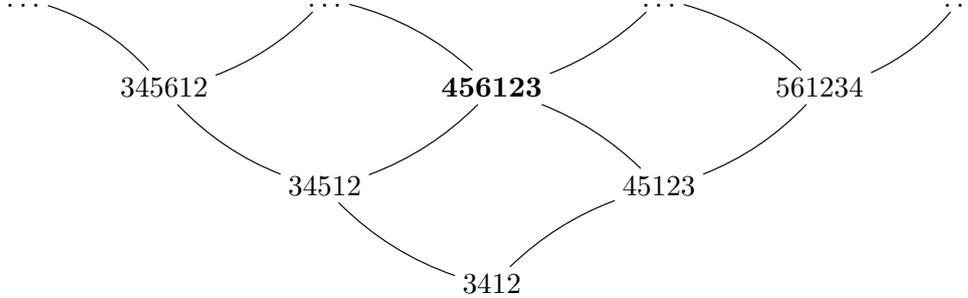
\begin{figure}[h!]
\centerline{
\xymatrix{
\dots & & \dots & & \dots & & \dots \\
 & 345612 \ar@/_/@{-}[dr] \ar@/_/@{-}[ul] \ar@/_/@{-}[ur] & & \mathbf{456123} \ar@/^/@{-}[dl] \ar@/^/@{-}[dr] \ar@/_/@{-}[ul] \ar@/_/@{-}[ur] & & 561234 \ar@/^/@{-}[dl] \ar@/_/@{-}[ul] \ar@/_/@{-}[ur] & \\
 & & 34512 \ar@/_/@{-}[dr] & & 45123 \ar@/_/@{-}[dl] & & \\
 & & & 3412 & & &
}
}
\caption{The set $\mathcal{A}_{2^+ 1^+}$.}\label{2+1+}
\end{figure}

\emph{Example.}\quad Consider the ball $\hat{B}_{1}^{(rd)}$,
whose clean compact peg basis is $\{1^- 2^-, 2^+ 1^{\bullet}, 2^{\bullet} 1^+ \}$.
Here we have
$\mathcal{M}_{1^- 2^-}= \left\lbrace 2143 \right\rbrace$,
$\mathcal{M}_{2^+ 1^{\bullet}}= \left\lbrace 231 \right\rbrace$ and
$\mathcal{M}_{2^{\bullet} 1^+}= \left\lbrace 312 \right\rbrace$,
so $B_1^{(rd)}=Av(2143,231,312)$ and the basis is indeed $\{ 2143,231,312\}$.
Notice that 2143 is the minimum of the set $\mathcal{A}_{1^- 2^-}$, and an analogous fact holds for 231 and 312 as well.
However, this is not true in general.
For instance, $2^+ 1^+$ belongs to the clean compact peg basis of $\hat{B}_{2}^{(rd)}$,
but the minimum permutation in $\mathcal{A}_{2^+ 1^+}$ is $3412$,
which has distance $2$ (unlike $2^+ 1^+$, which has distance 3).
So $3412 \notin \mathcal{M}_{2^+ 1^+}$.
Moreover (see Figure \ref{2+1+}) the permutations covering 3412 in $\mathcal{A}_{2^+ 1^+}$ are $34512$ and $45123$,
which have distance $2$ as well.
Therefore, in order to find the minimal permutations of $\mathcal{A}_{2^+ 1^+}$ at distance $3$,
we have to reach length $6$, where it can be shown that $456123$ is in fact minimal at distance $3$.
Thus $456123$ is a basis element corresponding to the clean compact peg basis permutation $2^+ 1^+$.

\bigskip

One nice feature of this approach is that the sets $\mathcal{A}_{\beta_i^{b_i}}$ are disjoint.
Thus each basis permutation is generated by precisely one of such sets.
On the other hand, it is still not clear if a single set $\mathcal{M}_i$ can lead to multiple basis permutations.

\bigskip

\textbf{Open problem 2.}\quad \emph{For a given clean compact peg basis permutation $\beta^b$,
determine $|\mathcal{M}_{\beta^b}|$. Is there always a single minimal permutation?}

\bigskip

\textbf{Open problem 3.}\quad \emph{Is the set $\mathcal{M}$ an antichain of permutations?
In other words, is $\mathcal{M}$ the standard basis?}

\bigskip

\textbf{Open problem 4.}\quad \emph{Determine a non trivial bound for the length of a basis permutation of $B_k ^{(rd)}(k)$,
for $k \ge 1$.}

\section{The prefix reversal model}

\subsection{Generating permutations}

In this section we focus on \textit{prefix reversals},
i.e. reversals that only affect the initial portion of a permutation.
Denoting with $prd$ the prefix reversal distance,
our first goal is to obtain a characterization of the balls $B_k ^{(prd)}$ in terms of generating permutations,
similarly to what we did in the case of the standard reversal distance.
As a first example, it is easy to see that $B_1 ^{(prd)}=Grid(1^-2^+)$,
so that $1^-2^+$ is the only generating permutation.
By specializing the recursive construction developed in the case of general reversals,
we can explicitly determine all the $(k+1)$-generating permutations starting from the $k$-generating ones.

\begin{prop}\label{prefixinfl}
Let $\pi^{\varepsilon}=\pi_1^{\varepsilon_1} \cdots \pi_n^{\varepsilon_n}$ be a generating permutation of $B_k ^{(prd)}$
and let $i \in \{1,\dots,n\}$.
\begin{enumerate}
\item If $\pi^{\varepsilon}=\alpha \ \pi_i^+ \ \beta$, for some peg permutations $\alpha,\beta$,
then $\pi_i^- \dot{\alpha}^R (\pi_i+1)^+ \dot{\beta}$ is a generating permutation of $B_{k+1}^{(prd)}$,
where $\dot{\alpha},\dot{\beta}$ are obtained from $\alpha,\beta$ (respectively) by
increasing by $1$ all the entries that are greater than $\pi_i$ and $\dot{\alpha}^R$ is the reversal of $\dot{\alpha}$.
\item If $\pi^{\varepsilon}=\alpha \ \pi_i^- \ \beta$, for some peg permutations $\alpha,\beta$,
then $(\pi_i+1)^+ \dot{\alpha}^R (\pi_i)^- \dot{\beta}$ is a generating permutation of $B_{k+1}^{(prd)}$,
where $\dot{\alpha},\dot{\beta}$ and $\dot{\alpha}^R$ are defined as above.
\end{enumerate}
Moreover, every $(k+1)$-generating permutation is obtained from a $k$-generating permutation by performing one of the two above constructions.
\end{prop}

\emph{Proof (sketch).} Since the prefix reversal model is a special case of the reversal one,
we have that every generating permutation for $B_{k+1}^{(prd)}$
can be obtained from a $k$-generating permutation $\pi^{\varepsilon}$ by choosing an index $i$,
then suitably inflating $\pi_i^{\varepsilon_i}$, according to its decoration,
and finally performing the prefix reversal operation that removes the new strip.
This is done in analogy with the construction described before Lemma \ref{prel},
with the only difference that the first index of the multiset has to be $1$
since we are considering prefix reversals.\cvd

The above proposition gives a recipe for constructing the generating permutations of $\hat{B}_{k+1}^{(prd)}$
starting from those of $\hat{B}_{k}^{(prd)}$.
Notice that, if $\pi^{\varepsilon}$ has length $m$,
then the permutations obtained with the previous construction have length $m+1$.
Since $B_1^{(prd)}=Grid(1^- 2^+)$,
a simple inductive argument shows that the generating permutations of $B_k^{(prd)}$ have length $k+1$.
Actually, we have something more,
which is the analogue of Corollary \ref{genblock} in the case of the prefix reversal model.
Since the proof is similar, we just state it.

\begin{theorem}
For every $k \geq 1$,
the generating set of $B_k ^{(prd)}$ is the set of all maximal clean compact peg permutations of $\hat{B}_{k}^{(prd)}$.
\end{theorem}

However, in the prefix case, we are able to count generating permutations.

\begin{theorem}
The generating set of $B_k ^{(prd)}$ has cardinality $k!$.
\end{theorem}

\emph{Proof.}\quad Observe that,
if $\pi^{\varepsilon}=\pi_1^{\varepsilon_1} \cdots \pi_{n}^{\varepsilon_{n}}$ is a $(k+1)$-generating permutation,
then it has been obtained from a $k$-generating permutation by one of the constructions
described in Proposition \ref{prefixinfl}.
However, $\pi^{\varepsilon}$ cannot be obtained in two different ways.
This can be shown by considering the elements $\pi_1^{\varepsilon_1}$:
\begin{enumerate}
\item if $\pi_1^{\varepsilon_1}=a^-$, for some $a \ge 1$,
then $\pi^{\varepsilon}$ is obtained as in $1.$ of Proposition \ref{prefixinfl}, when $\pi_i^{\epsilon_i}=a^+$;
\item if $\pi_1^{\varepsilon_1}=a^+$, for some $a \ge 1$,
then $\pi^{\varepsilon}$ is obtained as in $2.$ of Proposition \ref{prefixinfl}, when $\pi_i^{\epsilon_i}=(a-1)^-$.
\end{enumerate}

Since the above cases are disjoint,
we can conclude that $\pi^{\varepsilon}$ comes from a unique generating permutation of $B_k ^{(prd)}$
through the construction of Proposition \ref{prefixinfl}.
Thus, the total number of generating permutations of $B_{k+1}^{(prd)}$ is obtained by
multiplying the number of generating permutations of $B_k ^{(prd)}$ by the number of possible inflations of each of them,
which is $k+1$.
Since the generating set of $B_1 ^{(prd)}$ has cardinality $1$,
a simple inductive argument shows that the required cardinality is indeed $\prod_{i=1}^{k}{i}=k!$.\cvd

We have already observed that, for $k=1$, the generating set is $\{ 1^- 2^+ \}$.
For $k=2$, the generating set is $\{ 2^+ 1^- 3^+, 2^- 1^+ 3^+ \}$,
while for $k=3$ we obtain the set
$\{ 2^-3^+1^-4^+, 2^+3^-1^-4^+, 3^-1^+2^-4^+, 3^+2^-1^+4^+, 1^-3^+2^+4^+, 3^-1^-2^+4^+\}$.

\subsection{The bases}

Concerning the clean compact peg basis of $\hat{B}_{k}^{(prd)}$,
we are able to prove an analogue of Theorem \ref{basis}, although we need to take care of some exceptions.

\begin{lemma}\label{prf_breakp}
Let $\pi^{\varepsilon}$ be a clean compact peg permutation of length $n$ that does not end with its maximum value.
Then $prd(\pi^{\varepsilon})\geq n$.
\end{lemma}

\emph{Proof.} \quad Since any prefix reversal can change the number of breakpoints of a peg permutation by at most $1$
and a clean compact peg permutation of length $n$ has exactly $n-1$ breakpoints,
we have that $prd(\pi^{\varepsilon}) \ge n-1$.
Moreover, since $\pi^{\varepsilon}$ does not end with $n$,
in any sequence of prefix reversals that transforms the identity in $\pi^{\varepsilon}$
there must be a prefix reversal $R$ that moves $n$ from the end.
This means that $R$ reverses the entire permutation.
In particular, $R$ does not modify the number of breakpoints,
so we need at least $(n-1)+1=n$ reversals to obtain $\pi^{\varepsilon}$.\cvd

In the proof of the next theorem the software PermLab \cite{AL} has been of great help.

\begin{theorem}\label{exceptions}
Let $k \ge 0$ and $n=k+2$.
The following are clean compact peg basis permutations of $\hat{B}_{k}^{(prd)}$:
\begin{enumerate}
\item $$\begin{cases}
\Theta_e(n)= n^{\bullet} (n-2)^{\bullet} \dots 4^{\bullet} 2^{\bullet} 1^+ 3^{\bullet} \dots (n-3)^{\bullet} (n-1)^{\bullet}, \\
\Lambda_e(n)=(t+1)^+ t^{\bullet} (t+2)^{\bullet} (t-1)^{\bullet} \dots (n-1)^{\bullet} 2^{\bullet} n^{\bullet} 1^{\bullet},
\end{cases}$$
if $n$ is even and $t=\frac{n}{2}$;
\item $$\begin{cases}
\Theta_o(n)=  n^{\bullet} (n-2)^{\bullet} \dots 3^{\bullet} 1^- 2^{\bullet} 4^{\bullet} \dots (n-3)^{\bullet} (n-1)^{\bullet},\\
\Lambda_o(n)=t^- (t+1)^{\bullet} (t-1)^{\bullet} (t+2)^{\bullet} \dots (n-1)^{\bullet} 2^{\bullet} n^{\bullet} 1^{\bullet},
\end{cases}$$
if $n$ is odd and $t=\frac{n+1}{2}$.
\end{enumerate}

Moreover, each of the above peg permutations has distance $n=k+2$.

In the sequel, these peg permutations will be called \emph{exceptional}.
\end{theorem}

\emph{Proof.}\quad The proof will use induction on $k$.
If $k=0$, we get $\Theta_e(2)=2^{\bullet} 1^+$ and $\Lambda_e(2)=2^+ 1^{\bullet}$,
which are minimal clean compact peg permutations in the complement of $\hat{B}_{0}^{(prd)}$ and have distance $2$.
If $k=1$, we get $\Theta_o(3)=3^{\bullet} 1^- 2^{\bullet}$ and $\Lambda_o(3)=2^- 3^{\bullet} 1^{\bullet}$,
again minimal clean compact peg permutations in the complement of $\hat{B}_{1}^{(prd)}$ and having distance $3$.

If $k \ge 2$,
it is easy to observe that the above peg permutations are clean compact of length $n=k+2$
and do not end with their maximum value.
Thus, by Lemma \ref{prf_breakp}, they have distance at least $n=k+2$.

Next, we prove that the exceptional permutations are minimal in the complement of $\hat{B}_{k}^{(prd)}$.
In other words, we prove that each clean compact pattern of one of these permutations has distance at most $k$.
We first consider
$\Theta_e(n)=n^{\bullet} (n-2)^{\bullet} \dots 4^{\bullet} 2^{\bullet} 1^+ 3^{\bullet} \dots (n-3)^{\bullet}
(n-1)^{\bullet}$, with $n$ even.
There are several ways to obtain a clean compact pattern.
\begin{enumerate}

\item If we remove $i^{\bullet}$, with $i \neq 1,n$, in order to obtain a clean compact peg permutation,
we also need to remove either $i+1$ (if $i$ is even) or $i-1$ (if $i$ is odd);
otherwise, in fact, we would obtain the strip $(i-1),i$ or $i,(i-1)$, respectively.
We thus obtain the permutation
$(n-2)^{\bullet} (n-4)^{\bullet} \dots 4^{\bullet} 2^{\bullet} 1^+ 3^{\bullet} \dots (n-5)^{\bullet} (n-3)^{\bullet}=
\Theta_e(n-2)$, which has distance $n-2=k$ by induction hypothesis.

\item If we replace $1^+$ with $1^{\bullet}$,
we obtain the pattern
$n^{\bullet} (n-2)^{\bullet} \dots 4^{\bullet} 2^{\bullet} 1^{\bullet} 3^{\bullet}
\dots (n-3)^{\bullet} (n-1)^{\bullet}$,
which is not clean compact because of the strip $2^{\bullet} 1^{\bullet}$.
In order to obtain a clean compact pattern, we have to remove either $2^{\bullet}$ or $1^{\bullet}$, but in both cases,
after rescaling, we get the permutation
$(n-1)^{\bullet} (n-3)^{\bullet} \dots 3^{\bullet} 1^{\bullet} 2^{\bullet} \dots (n-4)^{\bullet} (n-2)^{\bullet}$,
which contains $1^{\bullet} 2^{\bullet}$.
Iterating this argument,
we get that the only clean compact pattern of $\Theta_e(n)$ obyained in this way is $1^{\bullet}\in \hat{B}_{k}^{(prd)}$.
The same happens if we choose to remove $1^+$.

\item Finally, if we remove the first element $n^{\bullet}$,
we obtain $(n-2)^{\bullet} \dots 4^{\bullet} 2^{\bullet} 1^+ 3^{\bullet} \dots (n-3)^{\bullet} (n-1)^{\bullet}$.
Notice that such a permutation ends with its maximum $(n-1)^{\bullet}$,
so it has the same distance as
$(n-2)^{\bullet} \dots 4^{\bullet} 2^{\bullet} 1^+ 3^{\bullet} \dots (n-3)^{\bullet}=\Theta_e(n-2)$
and we can conclude using the induction hypothesis.

\end{enumerate}

In a completely analogous way,
we can show that $\Lambda_e(n), \Theta_o(n)$ and $\Lambda_o(n)$
are minimal clean compact peg permutations in the complement of $\hat{B}_{k}^{(prd)}$ as well.

To conclude the proof, we now have to show that each of the exceptional permutations has distance exactly $n=k+2$.
We already know that such a distance is at least $n$, as a consequence of Lemma \ref{prf_breakp}.
Regarding the opposite inequality, consider for instance $\Theta_e(n)$.
Applying $2$ suitable prefix reversals to $\Theta_e(n)$,
we obtain the peg permutation $\Theta_e(n-2)$ followed by $(n-1)^{\bullet} n^{\bullet}$.
This is shown below:

\begin{equation}
\begin{split}
\Theta_e(n)=& \underbrace{n^{\bullet} (n-2)^{\bullet} \dots 4^{\bullet} 2^{\bullet} 1^+ 3^{\bullet} \dots
(n-3)^{\bullet} (n-1)^{\bullet}}_{PR_1} \\
\rightsquigarrow & \underbrace{(n-1)^{\bullet} (n-3)^{\bullet} \dots 3^{\bullet} 1^{-} 2^{\bullet} 4^{\bullet} \dots
(n-2)^{\bullet}}_{PR_2} n^{\bullet} \\
\rightsquigarrow &  (n-2)^{\bullet} \dots 4^{\bullet} 2^{\bullet} 1^+ 3^{\bullet} \dots
(n-3)^{\bullet} (n-1)^{\bullet} n^{\bullet}=\\
&=[\Theta_e(n-2)] \ (n-1)^{\bullet} n^{\bullet} := \hat{\Theta}.
\end{split}
\end{equation}

Notice that $prd(\hat{\Theta}) \le prd(\Theta_e(n-2))$, so, using the induction hypothesis:

$$prd(\Theta_e(n)) \le prd (\hat{\Theta})+2 \le prd(\Theta_e(n-2))+2=n-2+2=n,$$
as desired.
The same inequality can be proved in the same way for $\Theta_o(n)$.
Finally, consider the permutation $\Lambda_e(n)$.
Using a similar approach, we can apply a suitable prefix reversal as follows:

\begin{equation}
\begin{split}
\Lambda_e(n)=\underbrace{(t+1)^+}_{PR} t^{\bullet} (t+2)^{\bullet} (t-1)^{\bullet} \dots
(n-1)^{\bullet} 2^{\bullet} n^{\bullet} 1^{\bullet} \\
\rightsquigarrow (t+1)^- t^{\bullet} (t+2)^{\bullet} (t-1)^{\bullet} \dots
(n-1)^{\bullet} 2^{\bullet} n^{\bullet} 1^{\bullet} = \hat{\Lambda};
\end{split}
\end{equation}
the permutation $\hat{\Lambda}$ starts with the strip $(t+1)^- t^{\bullet}$, so, compacting and rescaling,
we obtain the peg permutation
$t^- (t+1)^{\bullet} (t-1)^{\bullet} (t+2)^{\bullet} \dots (n-2)^{\bullet} 2^{\bullet} (n-1)^{\bullet} 1^{\bullet}
=\Lambda_o(n-1)$.
Then, using again the induction hypothesis:

$$prd(\Lambda_e(n)) \le prd(\hat{\Lambda})+1 = prd(\Lambda_o(n-1))+1=(n-1)+1=n.$$

The case $\Lambda_o(n)$ can be dealt with in the same way.\cvd

\begin{cor}
\begin{itemize}
\item[(i)] If $n=k+2$ is even,
then $\Theta_e(n)$ and $\Lambda_e(n)$ are clean compact peg basis permutations for both
$\hat{B}_{k}^{(prd)}$ and $\hat{B}_{k+1}^{(prd)}$.
\item[(ii)] If $n=k+2$ is odd,
then $\Theta_o(n)$ and $\Lambda_o(n)$ are clean compact peg basis permutations for both
$\hat{B}_{k}^{(prd)}$ and $\hat{B}_{k+1}^{(prd)}$.
\end{itemize}
\end{cor}

\begin{lemma}\label{end_max}
Let $\pi^\varepsilon$ be a clean compact peg basis permutation of $\hat{B}_{k}^{(prd)}$ of length $k+2$.
Suppose $\beta^b$ is obtained from $\pi^{\varepsilon}$ by removing an element and suitably rescaling the remaining ones.
If $\beta^b$ is clean compact, then $\beta_{k+1}^{b_{k+1}} = (k+1)^{+ / \bullet}$.
\end{lemma}

\emph{Proof.} \quad Let $\pi^{\varepsilon} =\pi_1^{\varepsilon_1} \cdots \pi_{k+2}^{\varepsilon_{k+2}}$.
By the analogue of Lemma \ref{basis_prop} in the case of prefix reversal,
$\pi_{k+2}^{\varepsilon_{k+2}} \neq (k+2)^{+ / \bullet}$.
Moreover, by minimality of $\pi^{\varepsilon}$, it is $\beta^b \in \hat{B}_{k}^{(prd)}$,
so $\beta^b=\alpha^a \cdot [v]_{peg}$ for some $k$-generating permutation $\alpha^a$ and some inflation vector $v$.
Notice that $\alpha^a$ and $\beta^b$ have the same length $k+1$ and $\beta^b$ is clean compact,
hence $v_j=1$ for each $j=1,\ldots,k+1$.
Moreover, since it is easy to show that
the last element of any generating permutation of $\hat{B}_{k}^{(prd)}$ is its maximum decorated $+$,
we have that $\beta_{k+1}^{b_{k+1}}=(k+1)^{+ / \bullet}$, as desired.\cvd

\begin{theorem}
Every clean compact peg basis permutation of $\hat{B}_{k}^{(prd)}$ has length at most $k+2$.
\end{theorem}

\emph{Proof.} \quad The proof is identical to the first part of Theorem \ref{basis}.

\begin{theorem}\label{pref_bound}
Let $n=k+2$. Every clean compact peg basis permutation of $\hat{B}_{k}^{(prd)}$ has length at most $k+1=n-1$,
apart from the exceptional ones $\Theta_e(n), \Theta_o(n), \Lambda_e(n)$ and $\Lambda_o(n)$.
\end{theorem}

\emph{Proof.} \quad Suppose that
$\pi^{\varepsilon}=\pi_1^{\varepsilon_1} \cdots \pi_{n}^{\varepsilon_{n}}$ is a clean compact basis permutation of $\hat{B}_{k}^{(prd)}$ of length $n=k+2$.
By Lemma \ref{Reduced_pattern},
we know that it is possible to remove one element of $\pi^{\varepsilon}$ in such a way that the resulting peg permutation $\beta^{b}$ is clean compact
and thus, by the previous lemma, we have that $\beta_{n-1}^{b_{n-1}} = (n-1)^{+ / \bullet}$.
In particular, since $\pi_{n}^{\varepsilon_{n}} \neq n^{+ / \bullet}$,
$\beta^b$ can be obtained from $\pi^{\varepsilon}$ in one of the following (mutually exclusive) ways:
\begin{enumerate}
\item $\pi^{\varepsilon}$ ends with $n^{+ / - / \bullet} (n-1)^{+ / - / \bullet}$ and $\beta^b$ is obtained by removing one of the last two elements;
\item the last element of $\pi^{\varepsilon}$ is $(n-1)^{+ / \bullet}$ and $\beta^b$ is obtained by removing $n^{+/-/\bullet}$,
which is not the second to last element;
\item the second to last element of $\pi^\varepsilon$ is $n^{+ / \bullet}$ and $\beta^b$ is obtained by removing the last element, which is not $(n-1)^{+ / - / \bullet}$.
\end{enumerate}

Our goal is to show that we can almost always remove another element from $\pi^{\varepsilon}$ to obtain a clean compact peg permutation
that does not end with its maximum, which is a contradiction with Lemma \ref{end_max}.
The only exceptions are the exceptional permutations of Theorem \ref{exceptions}.

\begin{enumerate}

\item Suppose that $\pi^{\varepsilon}= \rho^p n^{+ / - / \bullet} (n-1)^{+ / - / \bullet}$, for some peg permutation $\rho^p$.
Notice that $\rho^p$ has to be clean compact, because $\pi^{\varepsilon}$ is clean compact and its last two elements are the largest elements of the permutation,
so no new strip can be created by removing any of them.
Then, again by Lemma \ref{Reduced_pattern}, we can remove an element of $\rho^p$ in order to obtain a clean compact peg permutation $\hat{\rho}^{\hat{p}}$.
Thus the peg permutation $\hat{\rho}^{\hat{p}} (n-1)^{+ / - / \bullet} (n-2)^{+ / - / \bullet}$ is clean compact,
has length $n-1$ and does not end with its maximum, as desired.

\item In the second case, we can remove the last element $(n-1)^{+ / \bullet}$ in most of the situations,
except when we have either $(n-2)^{+ / \bullet} n^{+ / \bullet}$ or $n^{- / \bullet} (n-2)^{- / \bullet}$.

If $\pi^{\varepsilon}$ contains $(n-2)^{+ / \bullet} n^{+ / \bullet}$, we can remove $(n-2)^{+ / \bullet}$ to get a clean compact peg permutation,
unless $(n-3)^{+ / \bullet}$ is immediately before $(n-1)^{+ / \bullet}$ in $\pi^{\varepsilon}$.
By iterating this argument, we find out that we are always able to remove an element other than $n^{+ / \bullet}$,
except when
$\pi^{\varepsilon}= 2^{+ / \bullet} 4^{+ / \bullet} 6^{+ / \bullet} \cdots n^{+ / \bullet} 1^{+ / \bullet} 3^{+ / \bullet} 5^{+ / \bullet}
\cdots (n-1)^{+ / \bullet}$, if $n$ is even, and when
$\pi^{\varepsilon}= 1^{+ / \bullet} 3^{+ / \bullet} 5^{+ / \bullet} \cdots n^{+ / \bullet} 2^{+ / \bullet} 4^{+ / \bullet} 6^{+ / \bullet}
\cdots (n-1)^+$, if $n$ is odd.
In such cases we can remove the element $1^{+ / \bullet}$ thus obtaining a clean compact peg permutation that does not end with its maximum, as desired.

Otherwise, if $\pi^{\varepsilon}$ contains $n^{- / \bullet} (n-2)^{- / \bullet}$, we can try to remove $(n-2)^{- / \bullet}$ and argue analogously,
obtaining the ``bad" peg permutations
$n^{- / \bullet} (n-2)^{- / \bullet} \dots 4^{- / \bullet} 2^{- / \bullet} 1^{+ / \bullet} 3^{+ / \bullet} \dots (n-3)^{+ / \bullet} (n-1)^{+ / \bullet}$,
if $n$ is even and $t=\frac{n}{2}$, and
$n^{- / \bullet} (n-2)^{- / \bullet} \dots 3^{- / \bullet} 1^{- / \bullet} 2^{+ / \bullet} 4^{+ / \bullet} \dots (n-3)^{+ / \bullet} (n-1)^{+ / \bullet}$,
if $n$ is odd and $t=\frac{n+1}{2}$.
This time it is impossible to find an element that can always be removed in order to find a clean compact peg permutation.
Notice that we have to decorate $1$ with $+$, if $n$ is even, and with $-$, if $n$ is odd.
If we decorate the other elements in the ``minimal" way, i.e. using $\bullet$,
we obtain, respectively, $\Theta_e(n)$ if $n$ is even and $\Theta_o(n)$ if $n$ is odd, as desired.

\item The third case is somehow symmetric to the second one.
Proceeding in an analogous way (just starting with the removal of the second to last element of $\pi^{\varepsilon}$),
we can find a clean compact peg permutation having the required properties,
except for the permutations
$(t+1)^{+ / \bullet} t^{- / \bullet} (t+2)^{+ / \bullet} (t-1)^{- / \bullet} \dots (n-1)^{+ / \bullet} 2^{- / \bullet} n^{+ / \bullet} 1^{- / \bullet}$,
if $n$ is even and $t=\frac{n}{2}$, and
$t^{- / \bullet} (t+1)^{+ / \bullet} (t-1)^{- / \bullet} (t+2)^{+ / \bullet} \dots (n-1)^{+ / \bullet} 2^{- / \bullet} n^{+ / \bullet} 1^{- / \bullet}$,
if $n$ is odd and $t=\frac{n+1}{2}$.
This time we have that $t+1$ has to be decorated with $+$, if $n$ is even, and $t$ with $-$, if $n$ is odd.
Decorating the remaining elements with $\bullet$, as above, we obtain $\Lambda_e(n)$ and $\Lambda_o(n)$, respectively.\cvd
\end{enumerate}

As an example, if we choose $k=3$, then we have that the exceptional permutations
$\Theta_o(5)=5^{\bullet} 3^{\bullet} 1^- 2^{\bullet} 4^{\bullet}$ and $\Lambda_o(5)=3^{-} 4^{\bullet} 2^{\bullet} 5^{\bullet}  1^{\bullet}$
are clean compact peg basis permutations for both $\hat{B}_{3}^{(prd)}$ and $\hat{B}_{4}^{(prd)}$;
in fact, they have distance $5$ and all their clean compact patterns have distance at most $3$.
Notice that also $\Theta_e(4)$ and $\Lambda_e(4)$ are clean compact peg basis permutations for $\hat{B}_{3}^{(prd)}$.

\bigskip

Regarding the (standard) basis of $B_k^{(prd)}$, all the results presented in the previous section can be adapted.
The main difference is that, in the case of the prefix reversal distance,
there are clean compact peg basis permutations of $\hat{B}_{k}^{(prd)}$ that have distance $k+2$,
so there can be (even large) gaps between their distance and the distance of their clean compact patterns.
For example, this is true for the exceptional permutations of Theorem \ref{exceptions}.
In fact, we can say something more precise about the basis elements that derive from the exceptional permutations.

\begin{theorem}
For a given peg permutation $\pi^{\varepsilon}$, consider the set:
$$\mathcal{M}_{\pi^{\varepsilon}}=\left\lbrace \gamma: \ \gamma \mbox{ is minimal such that } peg(\gamma)=\pi^{\varepsilon} \mbox{ and }
prd(\gamma) = prd(\pi^{\varepsilon}) \right\rbrace.$$

Also, denote with $min(\pi^{\varepsilon})$ the (standard) permutation obtained by inflating each element of $\pi^{\varepsilon}$ decorated $\bullet$ by $1$,
and the remaining ones by $2$. Then, for each $n$:
\begin{enumerate}
\item if $n$ is even,
$\mathcal{M}_{\Theta_e(n)} = \left\lbrace min(\Theta_e(n)) \right\rbrace$ and $\mathcal{M}_{\Lambda_e(n)} = \left\lbrace min(\Lambda_e(n)) \right\rbrace$;
in particular, $min(\Theta_e(n))$ and $min(\Lambda_e(n)$ are basis permutations for $B_{n-1}^{(prd)}$.
\item If $n$ is odd,
$\mathcal{M}_{\Theta_o(n)} = \left\lbrace min(\Theta_o(n)) \right\rbrace$ and $\mathcal{M}_{\Lambda_o(n)} = \left\lbrace min(\Lambda_o(n)) \right\rbrace$;
in particular, $min(\Theta_o(n))$ and $min(\Lambda_o(n))$ are basis permutations for $B_{n-1}^{(prd)}$.
\end{enumerate}
\end{theorem}

\emph{Proof.} (sketch)
First of all, notice that, in general,
$peg(min(\pi^{\varepsilon}))=\pi^{\varepsilon}$ and $min(\pi^{\varepsilon})$ is the minimal permutation with this property
(with respect to the pattern involvement relation).
We have $\gamma:=min(\Theta_e(n))=(n+1) \ (n-1)\dots 5 \ 3 \ 1 \ 2 \ 4 \ 6 \dots (n-2) \ n$.
We wish to show that $prd(\gamma)=prd(\Theta_e(n))=n$, which implies that $\mathcal{M}_{\Theta_e(n)} = \left\lbrace \gamma \right\rbrace$,
and that each pattern of $\gamma$ has distance at most $n-1$; this guarantees that $\gamma$ is a basis permutation.
We give the idea of the proof for $\Theta_e(n)$, the other cases being similar.

Define an \textit{adjacency} of $\gamma$ as a consecutive pair $\gamma_i \gamma_{i+1}$ in $\gamma$ such that $\gamma_{i+1}=\gamma_i+1$.
Then, by adapting Lemma \ref{end_max} to non-peg permutations, using adjacencies in place of breakpoints,
and observing that $\gamma$ does not end with its maximum, we can prove that $prd(\gamma) \ge n$.
Conversely, we can mimick the sorting procedure used in Theorem \ref{exceptions} for $\Theta_e(n)$
to obtain an analogous sorting sequence of length $n$ for $\gamma$, so $prd(\gamma) \le n$.
Therefore we have $prd(\gamma)=n$ and $\gamma$ is the minimal permutation such that $peg(\gamma)=\Theta_e(n)$,
so $\mathcal{M}_{\Theta_e(n)} = \left\lbrace \gamma \right\rbrace$, as desired.

Finally, in analogy with the proof of Theorem \ref{exceptions}, 
a simple case by case analysis shows that each pattern of $\gamma$ can be sorted using at most $n-1$ prefix reversals, 
meaning that $\gamma$ belongs to the basis of $B_{n-1}^{(prd)}$.\cvd

As a consequence of the previous theorem,
we have that the exceptional permutations of length $k+2$ are clean compact peg basis permutations of both $\hat{B}_{k}^{(prd)}$ and $\hat{B}_{k+1}^{(prd)}$,
but they contribute to the standard basis just for $\hat{B}_{k+1}^{(prd)}$.
For example, if $k$ is even, then we proved that $min(\Theta_o(k+1))$ and $min(\Lambda_o(k+1))$ are basis permutations for $B_k^{(prd)}$;
moreover, also $\Theta_e(k+2)$ and $\Lambda_e(k+2)$ are in the clean compact peg basis, but $min(\Theta_e(n-2))$ contains $min(\Theta_o(k+1))$ as a pattern,
so it cannot belong to the basis; analogously, $min(\Lambda_o(k+2) \ge min(\Lambda_o(k+1)$.
Consider for instance the ball $\hat{B}_{1}^{(prd)}$.
Theorem \ref{pref_bound} and a direct computation guarantee that its clean compact peg basis is the set $\{1^{\bullet} 2^-, 2^{\bullet}1^+, 2^+ 1^{\bullet} \}$
together with the exceptional permutations $\Theta_o(3)=3^{\bullet} 1^- 2^{\bullet}$ and $\Lambda_o(3)=2^- 3^{\bullet} 1^{\bullet}$.
Observe that $2^{\bullet}1^+ = \Theta_e(2)$ and $2^+ 1^{\bullet}=\Lambda_e(2)$ are the exceptional permutations of length $2$.
The corresponding sets of candidate basis permutations are
$\mathcal{M}_{1^{\bullet} 2^-}= \left\lbrace 132 \right\rbrace$, $\mathcal{M}_{2^+ 1^{\bullet}}= \left\lbrace 231 \right\rbrace$,
$\mathcal{M}_{2^{\bullet 1^+}}= \left\lbrace 312 \right\rbrace$, $\mathcal{M}_{3^{\bullet} 1^- 2^{\bullet}}=\{4213 \}$ and
$\mathcal{M}_{2^- 3^{\bullet} 1^{\bullet}}= \{ 3241 \}$;
since $min(\Theta_e(2))=312 \le 4213 = min(\Theta_o(3))$, $min(\Lambda_e(2))=231 \le 3241 = min(\Lambda_o(3))$,
the basis of $\hat{B}_{1}^{(prd)}$ is $\left\lbrace 132,231,312 \right\rbrace$.

Using the same approach,
but with considerably more efforts, a tedious computation shows that the basis of $B_{2}^{(prd)}$ is
$\left\lbrace 132,3241,3412,4213,4231 \right\rbrace$. Together with the results of \cite{HV},
this implies that $|Av(132,3241,3412,4213,4231)|=B_{2}^{(prd)}=n^2 +1$ (sequence A002522 in Sloane's Encyclopedia).

%


\end{document}